\documentclass[12pt]{article}
\usepackage{amsfonts, amsmath, latexsym, vmargin, epsfig}
\usepackage{epsf}
\usepackage{url}
\title{Cube packings, second moment and holes}
\setpapersize{custom}{21cm}{29.7cm}
\setmarginsrb{1.7cm}{3.5cm}{1.7cm}{3.2cm}{0pt}{0pt}{0pt}{5mm}
\def\QuotientSpace#1#2{\leavevmode\kern-.0em\raise.2ex\hbox{$#1$}\kern-.1em/\kern-.1em\lower.25ex\hbox{$#2$}}

\def\1t{1}
\def\2t{2}
\def\4t{4}

\newcommand { \ib }[1] {\textit{\textbf{#1}}}

\author{Mathieu DUTOUR\\
        \normalsize  Institut Rudjer Boskovi\'c, Zagreb\\
\and
        Yoshiaki ITOH\\
        \normalsize  Institute of Statistical Mathematics, Tokyo\\
\and
	Alexei POYARKOV\\
	\normalsize Moscow State University, Moscow\\
}

\begin{document}
\newcommand{\RR}{{\bf R}}
\newcommand{\NN}{{\bf N}}
\newcommand{\QQ}{{\bf Q}}
\newcommand{\CC}{{\bf C}}
\newcommand{\ZZ}{{\bf Z}}
\newcommand{\TT}{{\bf T}}
\newtheorem{proposition}{Proposition}
\newtheorem{theorem}{Theorem}
\newtheorem{corollary}{Corollary}
\newtheorem{lemma}{Lemma}
\newtheorem{problem}{Problem}
\newtheorem{conjecture}{Conjecture}
\newtheorem{claim}{Claim}
\newtheorem{remark}{Remark}
\newtheorem{definition}{Definition}
\newcommand{\qed}{\hfill $\Box$ }
\newcommand{\proof}{\noindent{\bf Proof.}\ \ }

\maketitle

\begin{abstract}
We consider tilings and packings of $\RR^d$ by integral translates of
cubes $[0,2[^d$, which are $4\ZZ^d$-periodic.
Such cube packings can be described by cliques of an associated graph, which
allow us to classify them in dimension $d\leq 4$.
For higher dimension, we use random methods for generating some examples.

Such a cube packing is called {\em non-extendible} if we cannot insert
a cube in the complement of the packing.
 In dimension $3$, there is a unique non-extendible cube
packing with $4$ cubes. We prove that $d$-dimensional cube packings with
more than $2^d-3$ cubes can be extended to cube tilings.
We also give a lower bound on the number $N$ of cubes of non-extendible
cube packings.

Given such a cube packing and $z\in \ZZ^d$, we denote by $N_z$
the number of cubes inside the $\4t$-cube $z+[0,4[^d$
and call {\em second moment} the average of $N_z^2$.
We prove that the regular tiling by cubes has maximal second moment and give
a lower bound on the second moment of a cube packing in terms of
its density and dimension.

\end{abstract}

\section{Introduction}

A {\em general cube tiling} is a tiling of $\RR^d$ by translates 
of the hypercube $[0,2[^d$, which we call a $2$-cube.
A {\em special cube tiling} is a tiling of $\RR^d$ by
integral translates of the hypercube $[0,2[^d$, which are $4\ZZ^d$-periodic.
An example of such a tiling is the {\em regular cube tiling} of
$\RR^d$ by cubes of the form $z+[0,2[^d$ with $z\in 2\ZZ^d$.

In dimension $1$, there is only one type of special cube tiling,
while in dimension $2$, two following types of special cube tilings exist:
\begin{center}
\begin{minipage}{5cm}
\centering
\resizebox{3.0cm}{!}{\includegraphics{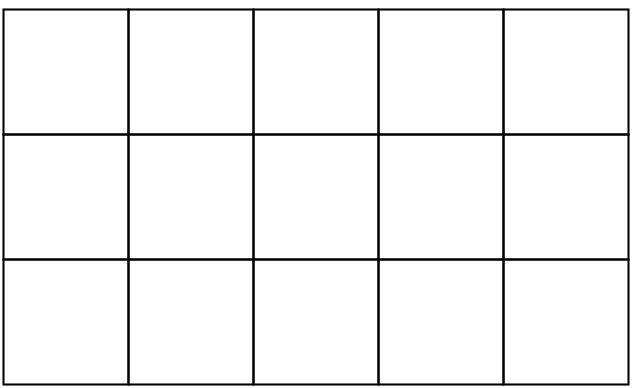}}\par
\end{minipage}
\begin{minipage}{5cm}
\centering
\resizebox{3.0cm}{!}{\includegraphics{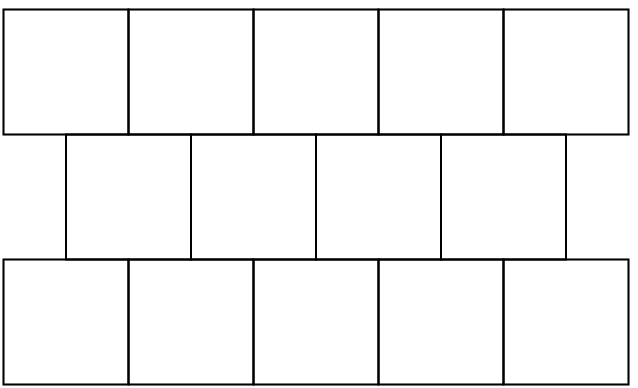}}\par
\end{minipage}

\end{center}
The Keller's cube tiling conjecture (see \cite{keller}) asserts
that any tiling of $\RR^d$ by translates of a unit cube admits
at least one face-to-face adjacency.
It is proved in \cite{szabo} that if this conjecture has a
counter example, then there is another counter example, which is
also a special cube tiling.
Using this, the Keller conjecture was solved
negatively for $d\geq 10$ in \cite{lagarias} and $d\geq 8$ in \cite{dim8}
(note that the conjecture is proved to be true for $d\leq 6$ in \cite{perron}).
Hence, special cube tilings, while seemingly limited objects have
a lot of combinatorial possibilities.
In the rest of this paper {\em cube tiling} stands for special cube
tilings and $N$ is the number of orbits of cubes
under the translation group $4\ZZ^d$. Another equivalent viewpoint
is to say that we are doing tilings of the torus
$\QuotientSpace{\RR^d}{4\ZZ^d}$ and $N$ is then the number of cubes
in this torus.

\begin{figure}
\begin{center}
\begin{minipage}{5.2cm}
\centering
\resizebox{5.0cm}{!}{\includegraphics[bb=105 260 494 559, clip]{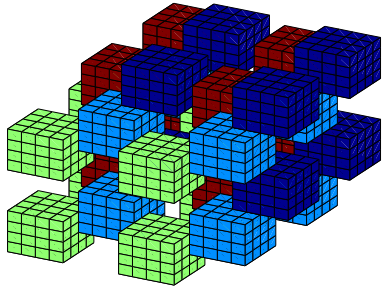}}\par
\end{minipage}
\end{center}
\caption{The unique non-extendible cube packing in dimension $3$}
\label{UniqueInex}
\end{figure}

A {\em cube packing} is a $4\ZZ^d$-periodic set of
integral translates of the $\2t$-cube, such that
any two cubes are non-intersecting.
In dimension $d\geq 3$, there exist cube packings, called 
{\em non-extendible}, which cannot be extended to a tiling
of the space (this first appear in \cite{lagariasEmail}).
In dimension $3$ this non-extendible packing is unique 
(see Figure \ref{UniqueInex}) and it is the source of much
of the inspiration of this paper.

In Section \ref{Algorithms}, following \cite{lagarias}, we present
a translation of the packing and tiling problems into clique problems
in graphs. Explicit methods, in GAP, are used up to $d=4$. For $d\geq 5$,
we use various random methods, in Fortran 90 and C++ for generating 
random cube packings.

Denote by $f(d)$ the smallest number of cubes, which form a non-extendible
cube packing. In Section \ref{Non-ExtendibleCubePackings}, we give some
lower and upper bounds on the value of $f(d)$.
In \cite{dip} it is proved that any cube packing of $[0,4[^d$ by cubes 
$[0,2[^d$ is extendible to a $4\ZZ^d$-periodic cube tiling of $\RR^d$.

If ${\cal CP}$ is a cube packing, denote by $hole({\cal CP})$ and call {\em hole}, its complement $\RR^d-{\cal CP}$.
We prove that if a cube packing has more than
$2^d-3$ cubes, then it is extendible to a tiling, i.e. that holes of volume
at most $3$ are fillable. We also obtain some conjecture on nonfillable
holes of volume at most $7$.

Given a cube packing ${\cal CP}$, the {\em counting function}
$N_z({\cal CP})$ is defined as the number of cubes of ${\cal CP}$
contained in $z+[0,4[^d$. We study its second moment 
in Section \ref{SecondMoment}. We prove, that the highest second
moment for tilings is attained for the regular cube tiling and give a
lower bound for the second moment of cube packings, in terms of its
dimension $d$ and number of cubes $N$.

\section{Algorithm for generating cube packings}\label{Algorithms}

Every $2$-cube of a $d$-dimensional cube packing is equivalent under
$4\ZZ^d$ to a cube with center in $\{0,1,2,3\}^d$.
Two $2$-cubes of centers $x$ and $x'$ do not
overlap if and only if there exist a coordinate $i$,
such that $|x_i-x'_i|=2$. So, one consider the graphs $G_d$ 
(introduced in \cite{lagarias}) with vertex-set $\{0,1,2,3\}^d$ and two
vertices being adjacent if and only if their associated cubes do
not overlap.
Cube packings correspond to cliques of $G_d$; they are non-extendible
if and only if the cliques are not included in larger cliques.
Cube tilings correspond to cliques of size $2^d$.

For a given $d$, the graph $G_d$ has a finite number of vertices and 
a symmetry group $Sym(G_d)$ of size $d! . 8^d$.
Hence, it is theoretically possible to do the enumeration of
the cliques of $G_d$.
The algorithm consists of using the set of all cliques with $N$ vertices,
considering all possibilities of extension, and then reducing by
isomorphism using $Sym(G_d)$ (the actual computation was done in GAP,
see \cite{Gap}).
The group $Sym(G_d)$ is presented as a permutation group in GAP and
the cliques as subsets of $\{1,\dots,4^d\}$. GAP uses backtrack search
for testing if two subsets are equivalent under $Sym(G_d)$ and
is hence, very efficient even for large values of $d$.
This enumeration is, in practice,  possible only for $d\leq 4$ due to
the huge number of cliques that appear.

For $d=2$, one finds only two
non-extendible cliques of $4$ vertices, i.e. two cube tilings.
For $d=3$, there is a unique non-extendible clique with $4$ vertices,
while there are $9$ orbits of non-extendible cliques with $8$ vertices
(i.e. cube tilings). 
For $d=4$, the computations are still possible and one finds the following
results with $N$ being the number of vertices of the clique.
\begin{equation*}
\begin{array}{|l|llllllllllllllll|}
\hline
N & 1& 2& 3& 4& 5& 6& 7& 8&  9& 10& 11& 12& 13& 14& 15& 16\\
\hline
\# \mbox{nonext.~orbit~cliques}
& 0& 0& 0& 0& 0& 0& 0& 38& 6& 24& 0&  71& 0&  0&  0&  744\\
\hline
\end{array}
\end{equation*}
Suppose that we have a cube tiling with two cube centers $x$ and $x'$ ,
satisfying to $x'=x+2e_i$ with $e_i$ being the $i$-th unit vector,
i.e. they have a face-to-face adjacency.
If one replace $x$, $x'$ by $x+e_i$, $x'+e_i$ and leave other centers
unchanged, then one obtains another cube tiling, which we call the
{\em flip} of the original cube tiling.
The enumeration strategy is then the following: take as initial list
of orbits the orbit of the regular cube tiling. For every orbit of
cube packing, compute all possible pairs $\{x,x'\}$, which allow to
create a new cube tiling. If the corresponding orbits of cube tilings
are new, then we insert them into the list of orbits.
Given a dimension $d$, consider the 
graph $Co_d$, whose vertex-set consists of all orbits of cube tilings
and put an edge between two orbits if one is obtained from the other by
a flipping. The above algorithm consists of computing the connected
component of the regular cube tiling in $Co_d$. Since the Keller conjecture is
false in dimension $d\geq 8$, we know that in those dimensions there are
some isolated vertices in the graph and so, the above algorithm does
not work. However, the graph $Co_d$ is connected for $d\leq 4$,
i.e. any two cube tilings in those dimensions can be obtained
by a sequence of flipping.
It is an interesting question to decide, in which dimension $d$ the graph
$Co_d$ is connected; the only remaining unsolved cases are $d=5,6,7$.

For dimensions $d\geq 5$, two above enumerative methods cannot work since there
are too much possibilities. Hence, we used random methods.
The {\em random packing} consist of selecting points, at
random, on $\{0,1,2,3\}^d$, so that the corresponding
$2$-cubes do not intersect, until one cannot do this any more.
This random packing algorithm creates non-extendible cube packings.

The actual algorithm for creating non-extendible cube packings is as follows:
the list $L$ of selected cubes is, initially, empty. One select at random
elements of $\{0,1,2,3\}^d$ and keep them if they are adjacent to preceding
elements of $L$. Of course not every trial works and as the space becomes
more and more filled, the number of random generation needed to get a
non-overlapping cube increase.
When this number has reached a certain level, we go to a second stage:
enumerate all possible extensions of the clique, that
we have, and work in this list by eliminating elements of it after choices are
made.
This algorithm has the advantage of enumerating the set $\{0,1,2,3\}^d$
only one time and is hence, relatively fast.

If one wants to find some packings with low density, then the above strategy
is not necessarily the best. The {\em greedy algorithm} consists of keeping
all $4^d$ elements in memory and at every step generate, say $20$ elements
and keep the one which cover the largest part of the remaining space.

Another possibility is what we call {\em Metropolis algorithm}
(see \cite{Liu}): we take
an non-extendible cube packing, remove a few cubes and rerun a random
generation from the remaining cubes. If obtained packing is better than
the preceding one, or not worse than a specified upper bound, then we keep
it; otherwise, we rerun the algorithm. This strategy allows to make a
random walk in the space of non-extendible cube packings and
is based on the assumption, that the best non-extendible cube packings
are not far from other, less good, non-extendible cube packings.

\section{Non-Extendible cube packings}\label{Non-ExtendibleCubePackings}

In dimension $1$ or $2$, any cube packing
is extendible to a cube tiling. 
The exhaustive enumeration methods of the preceding section show
that in dimension $3$, there is a unique non-extendible cube packing.
The set of its centers is, up to a symmetry of $G_3$:
\begin{equation*}
\{(0,0,0), (3,2,3), (2,1,1), (1,3,2)\}\;.
\end{equation*}
and its corresponding drawing is shown on Figure \ref{UniqueInex}.
Its space group symmetry is {\rm P4(1)32}, which is a chiral group.

We first concentrate on the problem of finding non-extendible cube packings
with the smallest number $f(d)$ of cubes.
From Section \ref{Algorithms}, we know that $f(1)=2$, $f(2)=4$, $f(3)=4$
and $f(4)=8$.


%
%
%



%


\begin{lemma}\label{MultiplicativeIneq}
For any $n,m\geq 1$, one has the inequality $f(n+m) \le f(n)f(m)$.
\end{lemma}
\proof Let $P_A$ and $P_B$ be non-extendible cube packings of $\RR^n$ and
$\RR^m$ with $f(n)$ and $f(m)$ cubes, respectively. Let
${\ib a}^k = (a_1^k, a_2^k, \ldots, a_n^k)$ and
${\ib b}^l = (b_1^l, b_2^l, \ldots, b_m^l)$ with
$1 \le k \le f(n)$, $1 \le l \le f(m)$ be the centers
of the $2$-cubes from $P_A$ and $P_B$.

Define $P$ to be the set of $2$-cubes $C^{kl}$ with centers 
${\ib c}^{kl} = (a_1^k, a_2^k, \ldots, a_n^k, b_1^l, b_2^l, \ldots, b_m^l)$
for $1 \le k \le f(n)$ and $1 \le l \le f(m)$. The size of $P$
is $f(n)f(m)$ and it is easy to check that $P$ is a packing.

Take a cube $D$ with center $\ib{d} = (d_1, d_2, \ldots, d_{n+m})$.
The vector $(d_1, \dots, d_n)$ overlaps with a $2$-cube, say $A^{k_0}$
in $P_A$, while the vector $(d_{n+1},\dots, d_{n+m})$ overlaps with
a $2$-cube, say $B^{l_0}$ in $P_B$. Clearly, $D$ overlaps with $C^{k_0l_0}$
and $P$ is non-extendible. \qed

Since, $f(3)=4$, one has $f(6) \le 16$.

\bigskip


A {\em blocking} set is a set $\{ {\ib v^j} \}$
of vectors in $\{0,1,2,3\}^d$,
such that for every other vector ${\ib v}$, there exist a $j$ such that the
$2$-cubes of center ${\ib v^j}$ and ${\ib v}$ overlap.
A priori, the $2$-cubes corresponding to the vector set $\{ {\ib v^j} \}$
can overlap; so, one has obviously $h(d)\leq f(d)$.

It is easy to see that $h(2)=3$ and that any blocking sets of size $3$
belong to one of two following orbits:
\begin{center}
\begin{minipage}{5cm}
\centering
\resizebox{3.0cm}{!}{\includegraphics{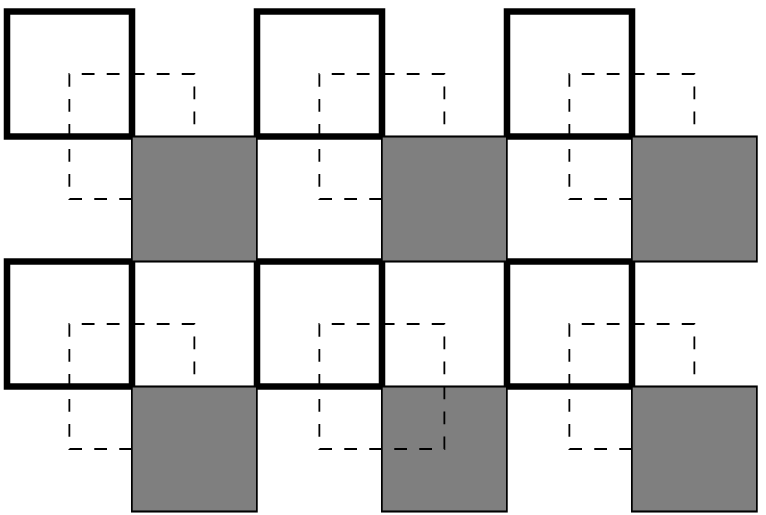}}\par
\end{minipage}
\begin{minipage}{5cm}
\centering
\resizebox{3.0cm}{!}{\includegraphics{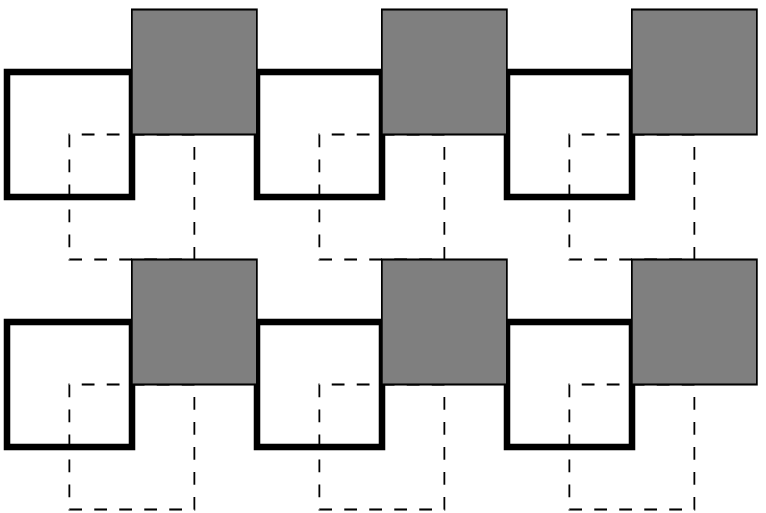}}\par
\end{minipage}

\end{center}
A slightly more complicated computation shows that $h(3)=4$
and that any blocking set of size $4$ belong to one of
three following orbits:
\begin{equation*}
\begin{array}{c}
\{(0,0,0), (1,1,1), (2,2,2), (3,3,3)\},\\
\{(0,0,0), (1,1,1), (2,2,3), (3,3,2)\},\\
\{(0,0,0), (3,2,3), (2,1,1), (1,3,2)\}
\end{array}
\end{equation*}

%
%


%
%
%
%
%
%

\begin{lemma}\label{KeyLemma}
Let $N$ satisfy the inequality
$\left\lfloor \frac{3N}{4} \right\rfloor < h(d)$, then one has $h(d+1) > N$.
\end{lemma}
\proof First $h(d)>N$ if and only if, for
any set $P$ of $N$ $\2t$-cubes, there exists a $\2t$-cube $D$, which
does not overlap with any $2$-cube from $P$.

Let $P$ be a  set of $N$ $\2t$-cubes in torus $T^{d+1}$.
Then at least $\left\lceil \frac{N}{4} \right\rceil$ centers of
them have $x_{d+1}=t$, for some $t\in \{0,1,2,3\}$.
Let us define another set $P'$ of vectors by removing those vectors
and the $d+1$-th coordinate for the remaining vectors.
Then $P'$ consists of at most 
$N - \left\lceil\frac{N}{4} \right\rceil = \left\lfloor \frac{3N}{4}
\right\rfloor$ $\2t$-cubes.
But $\left\lfloor\frac{3N}{4} \right\rfloor < h(d)$; so, there exists
a $2$-cube $C$ with center ${\ib c} = (c_1, c_2, \ldots, c_d)$,
which do not overlap with any $2$-cube in $P'$. But then
the $\2t$-cube with center $(c_1, c_2, \ldots, c_d, t+2)$ does not
overlap with any $2$-cube from $P$. \qed


\begin{theorem}\label{TheLowerBound}
For any $d\geq 1$, one has $h(d+1) \ge \left\lfloor \frac{4h(d)-1}{3}\right\rfloor +1.$
\end{theorem}
\proof Let $N = \left\lfloor \frac{4h(d)-1}{3} \right\rfloor$, then it holds:
$$
\left\lfloor \frac{3N}{4} \right\rfloor = \left\lfloor \frac{3
\left\lfloor \frac{4h(d)-1}{3} \right\rfloor}{4} \right\rfloor \le
\left\lfloor \frac{4h(d) -1}{4} \right\rfloor < h(d).
$$
And, from Lemma \ref{KeyLemma}, we have, that $h(d+1) > N$. \qed


Theorem \ref{TheLowerBound} does not allow to find an
asymptotically better lower bound on $f(d)$ than the trivial lower bound
$\lceil (\frac{4}{3})^d \rceil$. Note that using Lemma \ref{MultiplicativeIneq}
one proves easily that the limit $\beta=\lim_{d\to\infty}\frac{\ln\,f(d)}{d}$
exists. This limit satisfies to $\frac{4}{3}\leq e^{\beta}\leq \sqrt[3]{4}$.
The upper bound following from Lemma \ref{MultiplicativeIneq} and $f(3)=4$.
The determination of $\beta$ is open.

\begin{proposition}\label{ValueH4}
One has $h(4) = 7$.
\end{proposition}
\proof The following set of center coordinate proves that $h(4)\leq 7$.
\begin{equation*}
\begin{array}{c}
\{(0,0,0,0),
(1,1,1,1),
(2,2,2,2),
(3,3,3,3),
(0,0,1,1),
(1,1,2,2),
(2,2,3,3)\}
\end{array}
\end{equation*}
From Theorem \ref{TheLowerBound}, we have $h(4)\geq 6$.
Assume that $h(4)=6$ and take a blocking set of six $\2t$-cubes
with centers ${\bf a}^i = (a^i_1, a^i_2, a^i_3, a^i_4 )$, $1\le i \le 6$.

If three vectors ${\bf a}^{i_1}$, ${\bf a}^{i_2}$, ${\bf a}^{i_3}$
have equal coordinate $j$, then by a reasoning similar
to Lemma \ref{KeyLemma}, one finds a vector which does
not overlap with those vectors.

So, the above situation does not occur and for
every coordinate $j$, there exist two pairs
$\{{\bf a}^{i_1}, {\bf a}^{i_2}\}$, $\{{\bf a}^{i_3}, {\bf a}^{i_4}\}$,
which have equal $j$ coordinates.


We have two pairs $A$ and $B$ in first column.
Take a pair $A'$ in second column and assume that it does not
intersect with $A$. Denote by $P'$ the set of vectors obtained
by removing the vector corresponding to the sets $A$ and $A'$ and
the first and second coordinate of the remaining vectors. $P'$
is a set of two vectors in dimension $2$; hence, it is not blocking.
So, we can find a $2$-cube, which does not overlap with $P$.
So, any of six pairs from three other columns must intersect
with $A$ and $B$.

But we have only $4$ different
ways to intersect $A$ and $B$. So, two pairs
from column $2 - 4$ are equal. But, if two pairs are equal, then
they do not intersect, which is impossible. So, $h(4)>6$. \qed

Theorem \ref{TheLowerBound} and Proposition \ref{ValueH4} imply
the following inequalities:
\begin{equation*}
f(5)\geq h(5)\geq 10\mbox{~~~and~~~}f(6)\geq h(6)\geq 14
\end{equation*}

By running extensive random computation we found more than $140000$
non-extendible cube packings in dimension $5$ with $12$ cubes; they
belong to $203$ orbits.
Hence, it seems reasonable for us to conjecture that
in fact $f(5)=12$ and that the number of orbits of non-extendible
cube packings with $12$ cubes is ``small'', i.e. a few hundreds.

But dimension $6$ is already very different. We know that $f(6)\leq 16$
but we are unable to find by random methods a single non-extendible cube
packing with less than $20$ cubes.
\bigskip

%

We now consider cube packing with high density.

Take a cube packing of $\RR^d$ with center set $\{ {\bf x}^k\}$, 
$1\leq k\leq N$.
Select a coordinate $i$ and an index $j$ and form a cube packing
of $\RR^{d-1}$,
called {\em induced cube packing on layer $j$}, 
by selecting all ${\bf x}^k$ with $x^k_i=j,j+1\pmod 4$ and then
creating the vector $(x^k_1, \dots, x^k_{i-1},x^k_{i+1},\dots, x^k_n)$.

\begin{lemma}\label{Indution_Lemma}
If ${\cal CP}$ is a cube packing with $2^d-\delta$ cubes,
then its induced cube packings have at least $2^{d-1}-\delta$ cubes.
\end{lemma}
\proof Select a coordinate and denote by $n_j$ the number of $2$-cubes of ${\cal CP}$, with $x_i=j$.
One has $n_0+n_1+n_2+n_3=2^d-\delta$.

The number of $2$-cubes of the induced cube packing on layer $j$ is
$y_j=n_j+n_{j+1}$. One writes $y_j=2^{d-1}-\delta_j$ with $\delta_j\geq 0$,
since the induced cube packing is a packing. Clearly, one has
$\delta_0+\delta_1+\delta_2+\delta_3=2\delta$.

We have $n_j+n_{j+1}=2^{d-1}-\delta_j$; so, one gets, by subtracting
$n_j-n_{j+2}=\delta_{j+1}-\delta_j$, which implies:
\begin{equation*}
\delta_0-\delta_1+\delta_2-\delta_3=0\;.
\end{equation*}
Every vector $\Delta=(\delta_0,\delta_1, \delta_2,\delta_3)\in \ZZ_+^4$, satisfying
the above relation, can be expressed in the form 
$c_0(1,0,0,1)+c_1(1,1,0,0)+c_2(0,1,1,0)+c_3(0,0,1,1)$ with $c_j\in \ZZ_+$. This implies $\delta_j=c_j+c_{j+1}\leq \sum c_j=\delta$. \qed


\begin{theorem}
In dimension $d$, every cube packing with $2^d-\delta$ cubes for $\delta=1,2,3$ can be extended to a cube tiling.
\end{theorem}
\proof The proof is by induction on $d$. Take $d\geq 4$ and a cube packing
${\cal CP}$ with $2^d-\delta$ cubes and denote by $hole({\cal CP})$ its
hole in $\RR^d$. Let us consider the layering along the
coordinate $i$.
By Lemma \ref{Indution_Lemma}, the induced cube packings have
$2^{d-1}-\delta_j$ cubes with $\delta_j\leq 3$. So, one can complete them 
to form a cube packing of $\RR^{d-1}$. Denote by
${\cal CC}_i=[0,2[^{i-1}\times [0,1[\times [0,2[^{d-i}$ the half of a $2$-cube
cut along
the coordinate $i$. The induced cube packings are extendible by the induction
hypothesis. This
means that $hole({\cal CP})$ is the union of $2\delta$ cut cubes ${\cal CC}_i$.
Denote by $\Delta_i=(\delta_0,\delta_1, \delta_2,\delta_3)$
the corresponding vector; by the analysis of Lemma \ref{Indution_Lemma}
$\Delta_i=c_0(1,0,0,1)+c_1(1,1,0,0)+c_2(0,1,1,0)+c_3(0,0,1,1)$ for some
$c_i\in \ZZ_+$ with $\sum c_j=\delta$.

In the case $\delta=1$, it is clear that the only set, which for any $i$
can be written as ${\bf v}^{1,i}+{\cal CC}_i\cup {\bf v}^{2,i}+{\cal CC}_i$
for some vectors  ${\bf v}^{1,i}$, ${\bf v}^{2,i}$ is the $2$-cube itself.

If $\delta=2$, then, clearly, the vector $\Delta_i$ takes,
up to isomorphism, one of three different forms:
$(1,1,1,1)$, $(1,2,1,0)$ or $(2,2,0,0)$.

Suppose that for a given $i$, the vector $\Delta_i$ contains
the pattern $(0,1)$. This means that on one layer we have exactly one
translate, say ${\bf v}+{\cal CC}_i$, of ${\cal CC}_i$. Select any other
coordinate $i'$, ${\bf v}+{\cal CC}_i$ is splitted in two parts,
say ${\bf v}^1+{\cal CC}_{i'}$ and ${\bf v}^2+{\cal CC}_{i'}$
by the layers
along the coordinate $i'$. Since, an adjacent layer is completely filled,
this means that ${\bf v}^2={\bf v}^1\pm e_{i'}$.
Hence, they form a cube and the cube packing is extendible.

Suppose that for a given coordinate $i$, 
$\Delta_i=(x,x,0,0)$ with $x=2$ or $3$.
The $0$-th layer is filled with $x$
translates of set ${\cal CC}_i$.
Take another coordinate, say $i'$, and consider the
partition of $hole({\cal CP})$ into translates of ${\cal CC}_{i'}$.
By intersecting with the $0$-th layer, one obtains $2x$ intersections.
But since the third layer is full, it is necessary for the translate
of ${\cal CC}_{i'}$ to overlap
only on $1$-th layer. This means that they make a cube tiling.

The above considerations settle cases $(2,2,0,0)$ and $(1,2,1,0)$.
Now assume that for a given coordinate $i$, one has $\Delta_i=(1,1,1,1)$.
Assume also that the cube packing is non-extendible.
Take one translate ${\bf v}+{\cal CC}_i$ on layer $j$
in $hole({\cal CP})$. It is splitted in two parts by the translates of
${\cal CC}_{i'}$. Since we assume that the cube packing is non-extendible,
one of these translates overlaps on layer $j-1$ and the other one
on layer $j+1$.
One obtains a unique stair structure as illustrated below
in a two-dimensional section:
\begin{center}
\epsfig{height=30mm, file=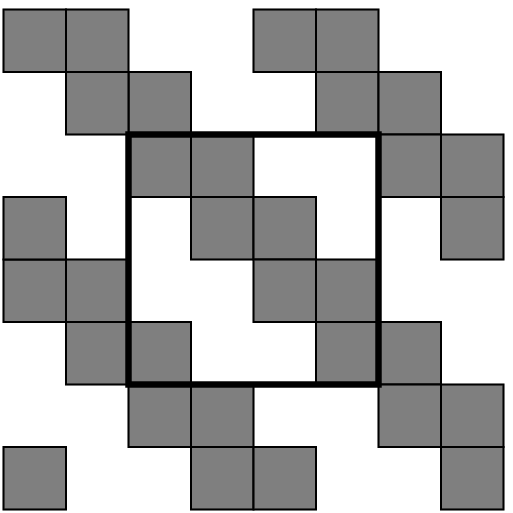}
\end{center}
Now select another coordinate $i''$ (since $d\geq 4$) and see that
$hole({\cal CP})$ cannot be decomposed into translates of ${\cal CC}_{i''}$. So, if
$\delta=2$, then all cube packings are extendible.

If $\delta=3$, then for a given coordinate $i$, one has clearly,
up to isomorphism, $\Delta_i$=$(3,3,0,0)$, $(2,1,1,2)$ or $(2,3,1,0)$.
The cases $(3,3,0,0)$ and $(2,3,1,0)$ are extendible by the above analysis.
Let us consider the case
$(2,1,1,2)$ and assume that the cube packing is non-extendible.
The $1$-th and $2$-th layers consist of translates of ${\cal CC}_i$, which
we write as ${\bf v}^1+{\cal CC}_i$ and ${\bf v}^2+{\cal CC}_i$.
The translate ${\bf v}^2+{\cal CC}_i$ is
splitted in two by the translate of ${\cal CC}_{i'}$ appearing in the decomposition
of $hole({\cal CP})$ along coordinate $i'$. If those translates spilled only on
the $0$-th layer or $2$-th layer, then one has a cube, which is excluded.
So, they spill on $0$-th and $2$-th layers. This implies that
${\bf v}^2={\bf v}^1\pm e_{i'}$. But this is impossible, since $i'$ is arbitrary. So, the cube packing is extendible. \qed



Given a $d$-dimensional non-extendible cube packing with $2^d-\delta$,
its {\em lifting} is a $d+1$-dimensional non-extendible cube packing
obtained by adding a layer of cube tiling; the iteration of lifting is
also called lifting.
\begin{conjecture}
Take ${\cal CP}$ a non-extendible cube packing with $2^d-\delta$
cubes. On its hole we conjecture:
\begin{enumerate}
\item If $\delta=4$ then $hole({\cal CP})$ is obtained as the hole of the
lifting of the unique non-extendible cube packing in dimension $3$.
\item The case $\delta=5$ does not occur.
\item If $\delta=6$, then $hole({\cal CP})$ is obtained as the hole of
the lifting of one of two non-extendible cube packing in dimension $4$.

\item If $\delta=7$, then $hole({\cal CP})$ is obtained as the hole of
the lifting of a non-extendible cube packing in dimension $4$.
\end{enumerate}
\end{conjecture}
This conjecture is supported by extensive numerical computations.
We can obtain an infinity of non-extendible cube packings with
$2^d-8$ cubes by doing layering of two $(d-1)$-dimensional
non-extendible cube packings with $2^{d-1}-4$ cubes.
This phenomenon does not appear for non-extendible cube packings
with $2^d-9$ cubes, but we are not able to state a reasonable conjecture
for this case.




\section{The second moment}\label{SecondMoment}

Given a cube packing ${\cal CP}$ and $z\in \ZZ^d$,
$N_z({\cal CP})$ is defined as the number of $2$-cubes
of ${\cal CP}$ contained in $z+[0,4[^d$.

Given a $4\ZZ^d$-periodic function $f$, its average is
\begin{equation*}
E(f)=\frac{1}{4^d}\sum_{z\in \{0,1,2,3\}^d} f(z)\;.
\end{equation*}
We denote $m_i({\cal CP})$ the $i$-moment of ${\cal CP}$, i.e.
the average of  $N_z^i({\cal CP})$.



\begin{theorem}
Let ${\cal CP}$ be a cube packing with $N$ cubes. One has:
\begin{equation*}
m_1({\cal CP})=(\frac{3^d}{4^d})N\mbox{~~~and~~~}
m_1({\cal CP})+N(N-1)2^{-d}+2^{-d} d\{2q(q-1)+rq\}\leq m_2({\cal CP})
\end{equation*}
with $N=4q+r$, $0\leq r\leq 3$.

\end{theorem}
\proof Take $N$ $2$-cubes $A^1$, \dots, $A^N$ with centers ${\bf a}^1$, \dots, ${\bf a}^N$.
The $\4t$-cube ${\ib a}+[0,4[^d$ with corner $(a_1,\dots, a_d)$
contains the $2$-cube with center ${\bf b}=(b_1,\dots, b_d)$ if and only
if $a_i\not= b_i$ for every $i$.
Take all $\4t$-cubes $C_1$, \dots, $C_{4^d}$.

Every $2$-cube $A^i$ is contained in $3^d$ $4$-cubes $C_k$.
Denote by $n_j$ the number
of $2$-cubes $A^i$, contained in the $4$-cube $C_j$.
By definition, the first moment has the expression:
\begin{equation*}
m_1({\cal CP})=\frac{1}{4^d}\sum_k n_k=\frac{1}{4^d}(3^dN)\;.
\end{equation*}
The second moment is equal to $m_2({\cal CP})=\frac{1}{4^d}\sum_k n^2_k$.
Let $t_{ij}$ be the numbers of $\4t$-cubes containing the $2$-cubes $A^i$
and $A^j$. One has the relation:
\begin{equation*}
\sum_{1\leq i < j\leq N}t_{ij}=\sum_{k=1}^{4^d} \frac{n_k(n_k-1)}{2},
\end{equation*}
which implies $4^d m_1({\cal CP})+2\sum t_{ij}=4^d m_2({\cal CP})$.
Let us denote by $\mu_{ij}$ the number of equal coordinates of the
centers ${\ib a}^i$ and ${\ib a}^j$.
Then one has
\begin{equation*}
t_{ij}=(\frac{3}{2})^{\mu_{ij}} 2^d\geq 2^d+2^{d-1}\mu_{ij}
\end{equation*}
The above inequality becomes an equality for $\mu_{ij}=0$ or $1$.
Summing over $i$ and $j$ one obtains
\begin{equation*}
\sum_{1\leq i<j\leq N} t_{ij}\geq N(N-1)2^{d-1}+2^{d-1}\sum_{1\leq i<j\leq N} \mu_{ij}
\end{equation*}
Let us denote by $R_l$ the number of equal pairs in column $l$.
By definition, one has clearly:
\begin{equation*}
\sum_{1\leq i<j\leq N} \mu_{ij}=\sum_{l=1}^d R_l \;.
\end{equation*}
Let us fix a coordinate $l$ and denote by $d_u$ the number of entries equal to $u$ in column $l$. One has, obviously:
\begin{equation*}
R_l=\sum_{u=0}^3 \frac{d_u(d_u-1)}{2},\;\;d_u\geq 0\mbox{~~and~~}\sum_{u=0}^3 d_u=N\;.
\end{equation*}
The Euclidean division $N=4q+r$ and elementary optimization, with respect to the constraints, allow us to write:
\begin{equation*}
R_l\geq 2q(q-1)+rq
\end{equation*}
The proof follows by combining all above elements. \qed

Note that the value of $m_1({\cal CP})$ was already obtained in \cite{dip}.
For a fix $d$ and $N$, we do not know which cube packing minimize the
second moment. However, we can characterize which cube tilings have
the highest second moment in Theorem \ref{MaxSecondMoment}.

%

Consider the following space of functions:
\begin{equation*}
{\cal G}=\left\lbrace\begin{array}{c}
f:\{0,1,2,3\}^d\rightarrow \RR.\\
\forall x\in\{0,1,2,3\}^d\mbox{~one~has~}\sum_{x+\{0,1\}^d}f(x)=1\mbox{~and~}f(x)\geq 0
\end{array}\right\rbrace\;.
\end{equation*}
It is easy to see that cube tilings correspond to $(0,1)$ vector in ${\cal G}$. Therefore, the problem of minimizing the second moment over cube tilings is an integer programming problem for a convex functional.


\begin{theorem}\label{MaxSecondMoment}
The regular cubic tiling is the cube tiling with highest second moment.
\end{theorem}
\proof Given a function $f\in {\cal G}$, let us define 
\begin{equation*}
M_i(f)(x)=\left\lbrace\begin{array}{rcl}
f(x)+f(x+e_i)&\mbox{~if~}& x_i=0\mbox{~or~}2\\
0            &\mbox{~if~}& x_i=1\mbox{~or~}3\;.
\end{array}\right.
\end{equation*}
The function $M_i(f)$ belongs to ${\cal G}$. Geometrically $M_i(f)$ is the cube packing obtained by merging two induced cube packing on coordinate $i$ and layer $0$ and $2$.
We will prove $E(N_z(M_i(f))^2)\geq E(N_z(f)^2)$.
Without loss of generality, one can assume, $i=1$.

The key inequality, used in computation below, is:
\begin{equation*}
\begin{array}{c}
(x_0+x_1+x_2)^2+(x_1+x_2+x_3)^2+(x_2+x_3+x_0)^2+(x_3+x_0+x_1)^2\\
\leq 2(x_0+x_1+x_2+x_3)^2+(x_0+x_1)^2+(x_2+x_3)^2\mbox{~if~}x_i\geq 0.
\end{array}
\end{equation*}
Define $f_{z_2}(z_1)=\sum_{u_2\in \{0,1,2\}^{d-1}}f(z_1, z_2+u_2)$ and obtain:
\begin{equation*}
\begin{array}{rcl}
4^d E(N_z(M_1(f))^2)
&=&\sum_{z \in \{0,1,2,3\}^d} (\sum_{u\in \{0,1,2\}^d} M_1(f)(z+u))^2\\
&=&\sum_{z_1=0}^3 \sum_{z_2\in \{0,1,2,3\}^{d-1}} (\sum_{u_1=0}^2 \sum_{u_2\in \{0,1,2\}^{d-1}} M_1(f)(z_1+u_1, z_2+u_2))^2\\
&=&\sum_{z_2\in \{0,1,2,3\}^{d-1}}\sum_{z_1=0}^3 (\sum_{u_1=0}^2 M_1(f_{z_2})(z_1+u_1))^2\\
&=&\sum_{z_2\in \{0,1,2,3\}^{d-1}}\{ 2(\sum_{u_1=0}^3 f_{z_2}(u_1))^2+(\sum_{u_1=0}^1 f_{z_2}(u_1))^2+(\sum_{u_1=2}^3 f_{z_2}(u_1))^2\}\\
&\geq&\sum_{z_2\in \{0,1,2,3\}^{d-1}}\sum_{z_1=0}^3(\sum_{u_1=0}^2 f_{z_2}(z_1+u_1))^2=4^d E(N_z(f)^2)
\end{array}
\end{equation*}
Hence, using the operation $M_1 \dots M_d$, we can only increase the
second moment. So, one gets:
\begin{equation*}
E(N_z(M_1\dots M_d(f))^2)\geq E(N_z(f)^2)\mbox{~for~all~}f\in {\cal G}\;.
\end{equation*}
It is easy to see that $M_1\dots M_d(f)$ is the function with
$f(x)=1$ if $x$ is a $(0,2)$ vector and $0$, otherwise; hence,
it corresponds to a regular cube tiling. \qed

Note that it is easy to see that $m_2=(\frac{5}{2})^d$ for the regular cube tiling.


\end{document}